\documentclass[12pt,twoside,reqno]{amsart}
\usepackage{graphicx}
\usepackage{amsmath}
\usepackage{amsfonts}
\usepackage{amssymb}
\usepackage{amscd}
\usepackage{amsmath}
\usepackage{amscd}
\usepackage{amsthm}
\usepackage{amssymb}
\usepackage{amsfonts}
\usepackage{epsfig}
\usepackage{amsmath,amscd}
\usepackage{graphicx}
\usepackage{geometry,enumerate,latexsym,tabularx,shapepar}
\usepackage{pslatex}
\usepackage[all,2cell,knot,poly]{xy} \UseAllTwocells \SilentMatrices
\usepackage{slashbox}
\input{xypic}
\xyoption{all}
\numberwithin{equation}{section}
\newcommand{\ie}{{\it i.e.\/}\ }

\newtheorem{thm}{Theorem}[section]

\newtheorem{prop}[thm]{Proposition}
\theoremstyle{definition}
\newtheorem{defn}[thm]{Definition}
\theoremstyle{remark}

\numberwithin{equation}{section}
\newtheorem{exa}[thm]{Example}

\newcommand{\bq}{\begin{eqnarray}}
\newcommand{\nq}{\end{eqnarray}}

\newcommand{\R}{\mathbb{R}}
\newcommand{\s}{\mathbf{S^3}}

\def\R{{\mathbb R}}

\begin{document}

\title[Khovanov-Kauffman Homology for embedded Graphs ]{Khovanov-Kauffman Homology for embedded Graphs }%
\author{Ahmad Zainy Al-Yasry} %

\address{University of Baghdad}%
\email{ahmad@azainy.com}%


\date{\today}%


\begin{abstract}
 A discussion given to the question of extending Khovanov homology from links to embedded graphs, by using the
 Kauffman topological invariant of embedded graphs by associating family of links and
 knots  to a such graph by using some local replacements at each vertex in the graph.
 This new concept of Khovanov-Kauffman homology of an embedded graph constructed to be the sum of the
Khovanov homologies of all the links and knots associated to this graph.
\end{abstract}

\maketitle

\section{Introduction}

The idea of categorification the Jones  polynomial is known by
Khovanov Homology for links which is a new link invariant introduced by
Khovanov \cite{kh}, \cite{Ba}. For each link $L$ in $S^3$ he defined a graded chain complex,
with grading preserving differentials, whose graded Euler characteristic is
equal to the Jones polynomial of the link $L$. The idea of Khovanov Homology
for graphs arises from the same idea of Khovanov homology for links by
the categorifications the chromatic polynomial of graphs.
This was done by L. Helme-Guizon and Y. Rong \cite{laur},
for each graph G, they defined a graded chain complex whose graded
Euler characteristic is equal to the chromatic polynomial of G.
In our work we try to recall, the Khovanov homology for links.\\
We discuss the question of extending Khovanov homology from links to
embedded graphs. This is based on a result of Kauffman that constructs a topological
invariant of embedded graphs in the 3-sphere by associating to such a graph
a family of links and knots obtained using some local replacements at
each vertex in the graph. He showed that it is a topological invariant by showing
that the resulting knot and link types in the family thus constructed are invariant
under a set of Reidemeister moves for embedded graphs that determine the
ambient isotopy class of the embedded graphs. We build on this idea and simply
define the Khovanov homology of an embedded graph to be the sum of the
Khovanov homologies of all the links and knots in the Kauffman invariant
associated to this graph. Since this family of links and knots is a topologically invariant,
so is the Khovanov-Kauffman homology of embedded graphs defined in this manner. We close this paper
by giving an example of computation of Khovanov-Kauffman homology for an embedded graph
using this definition.\\
\\
\textbf{\textrm{Acknowledgements:}} The author would like to express his deeply grateful to Prof.Matilde Marcolli
for her advices and numerous fruitful discussions. He is also thankful to Louis H. Kauffman and Mikhail Khovanov
for their support words and advices. Some parts of this paper done in Max-Planck Institut f\"ur Mathematik (MPIM), Bonn, Germany,
the author is kindly would like to thank MPIM for their hosting and subsidy during his study there.

\section{Khovanov Homology}

In the following we recall a homology theory for knots and links embedded in the
3-sphere. We discuss later  how to extend it to the case of embedded
graphs.
\\
\subsection{Khovanov Homology for links}
In recent years, many papers have appeared that discuss properties of Khovanov Homology theory,
which  was introduced in \cite{kh}. For each link  $ L \in \s $, Khovanov constructed
a bi-graded chain complex associated with the diagram $D$ for this link
$L$ and applied homology to get a group
$Kh^{i,j}(L)$, whose Euler characteristic is the normalized Jones polynomial.
$$ \sum_{i,j}(-1)^i q^j dim(Kh^{i,j}(L))=J(L)$$
He also proved that, given two diagrams $D$ and $D'$ for the same link, the corresponding
chain complexes are chain equivalent, hence, their homology groups are
isomorphic. Thus, Khovanov homology is a link invariant.

\subsection{The Link Cube}\label{licube}
\begin{figure}[htp]
\begin{center}
\includegraphics[width=6.5cm]{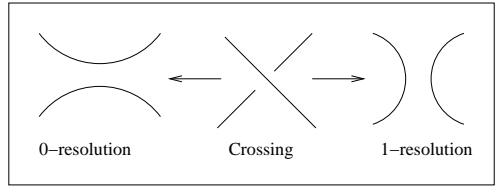}
\caption{0 and 1- resolutions to each crossing }\label{fig11}
\end{center}
\end{figure}

Let $L$ be a link with $n$ crossings. At any small neighborhood of a crossing we can
replace the crossing by a pair of parallel arcs and this operation is called a resolution.
There are two types of these resolutions called $0$-resolution (Horizontal resolution) and
$1$-resolution (Vertical resolution) as illustrated in figure (\ref{fig11}).\\
\\
We can construct a
$n$-dimensional cube by applying the $0$ and $1$-resolutions $n$ times to each crossing to
get $2^n$ pictures called smoothings (which are one dimensional manifolds) $S_\alpha$.
Each  of these can be indexed by a word $\mathbf{\alpha}$ of $n$ zeros and ones, \ie
$\mathbf{\alpha}  \in \{0,1\}^n$.
Let $\xi$ be an edge of the cube between two smoothings $S_{\alpha_{1}}$ and $S_{\alpha_{2}}$,
where $S_{\alpha_{1}}$ and $S_{\alpha_{2}}$ are identical smoothings except for a small
neighborhood around the crossing that changes from $0$ to $1$-resolution. To each edge
$\xi$ we can assign a cobordism $\Sigma_{\xi}$ (orientable surface whose boundary is the
union of the circles in the smoothing at either end)$$\Sigma_{\xi}: S_{\alpha_{1}} \longrightarrow S_{\alpha_{2}}$$
This $\Sigma_{\xi}$ is a product cobordism except in the neighborhood of the crossing, where it is  the obvious
saddle cobordism between the $0$ and $1$-resolutions.
Khovanov constructed a complex by applying a $1+1$-dimensional TQFT (Topological Quantum Field Theory)
which is a monoidal functor,
by replacing each vertex $S_\alpha$  by a graded vector space  $V_{\alpha}$
and each edge (cobordism) $\Sigma_{\xi}$ by a linear map $d_{\xi}$, and we set the group $CKh(D)$ to be the direct sum
of the graded vector spaces for all the vertices and the differential on the summand
$CKh(D)$ is a sum of the maps  $d_{\xi}$ for all edges $\xi$ such that Tail($\xi$)$=\alpha $ \ie
\begin{equation}\label{dv}
d^i(v)=\sum_{\xi} sign(-1)d_\xi(v)
\end{equation}
where $v \in V_\alpha \subseteq CKh(D)$ and $sign (-1)$  is chosen such that $d^2=0$.\\
An element of $CKh^{i,j}(D)$ is said to have homological grading $i$ and $q$-grading $j$ where
\begin{equation}\label{76}
i=|\alpha|-n_-
\end{equation}
\begin{equation}\label{77}
j=deg(v)+i+n_-+n_+
\end{equation}
for all $v \in V_\alpha \subseteq CKh^{i,j}(D)$, $|\alpha|$ is the number of 1's in $\alpha$, and
$n_-$, $n_+$ represent the number of negative and positive crossings respectively in the diagram $D$.

\subsection{Properties}\cite{pt1}, \cite{lee}
Here we give some properties of Khovanov homology.
\begin{prop}\label{1234}
\begin{enumerate}
\item If $D'$ is a diagram obtained from $D$ by the application of a Reidemeister moves then
the complexes $ (CKh^{*,*}(D))$ and $ (CKh^{*,*}(D'))$ are homotopy equivalent.
\item For an oriented link $ L$ with diagram D, the graded Euler characteristic satisfies
\begin{equation}\label{12}
\sum (-1)^i qdim(Kh^{i,*}(L))=J(L)
\end{equation}
where $ J(L)$
is  the normalized Jones Polynomials for a link $L$ and $$\sum (-1)^i qdim(Kh^{i,*}(D))=\sum(-1)^i qdim(CKh^{i,*}(D))$$
\item Let $L_{odd}$ and $L_{even}$  be two links with odd and even number of components then $Kh^{*,even}(L_{odd})=0$ and
$Kh^{*,odd}(L_{even})=0$\\
\item
For two oriented link diagrams $D$ and $D'$, the chain complex of the disjoint union $ D \sqcup D'$ is given by
\begin{equation}\label{123}
CKh(D \sqcup D') =CKh(D) \otimes CKh (D').
\end{equation}
\item For two oriented links  $L$ and $L'$, the Khovanov homology of the disjoint union $ L \sqcup L'$ satisfies
$$Kh(L \sqcup L') = Kh(L) \otimes Kh (L').$$
\item Let $D$ be an oriented link diagram of a link $L$ with mirror image $ D^m$ diagram of the mirror link $ L^m$.
Then the chain complex $CKh(D^m)$ is isomorphic to the dual of $CKh(D)$ and $$ Kh(L) \cong Kh(L^m)$$
\end{enumerate}
\end{prop}
\section{Khovanov-Kauffman Homology for Embedded Graphs (KKh(G)) }

\subsection{Kauffman's invariant of Graphs}
\begin{figure}
\begin{center}
\includegraphics[width=6cm]{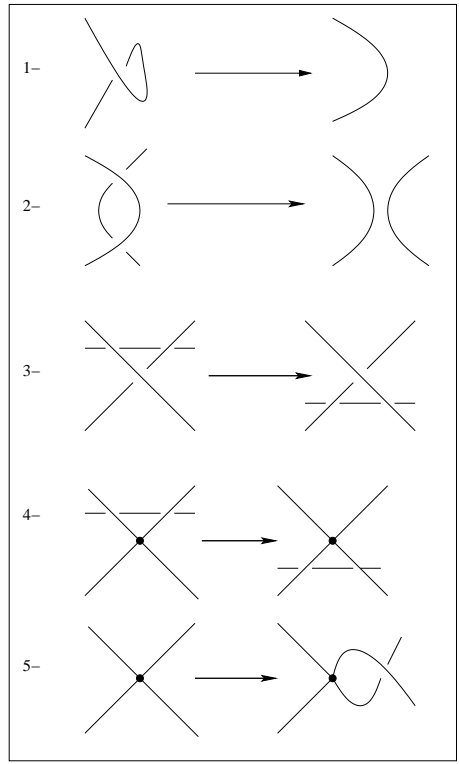}
\caption{Generalized Reidemeister moves by Kauffman}\label{fig21}
\end{center}
\end{figure}

\begin{figure}
\begin{center}
\includegraphics[width=8cm]{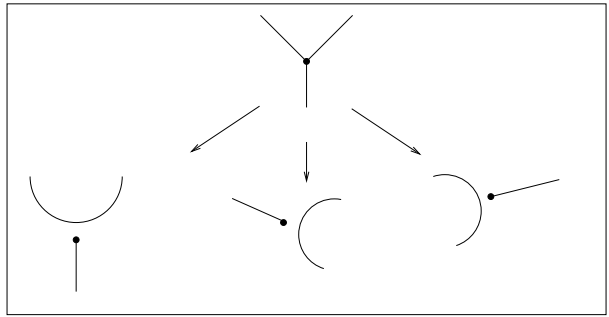}
\caption{local replacement to a vertex in the graph G}\label{fig18}
\end{center}
\end{figure}

\begin{figure}
\begin{center}
\includegraphics[width=8cm]{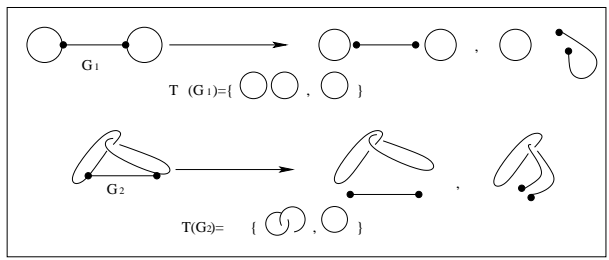}
\caption{Family of links associated to a graph}\label{fig19}
\end{center}
\end{figure}

We give now a survey of the Kauffman theory and show how to associate to an embedded
graph in $\s$ a family of knots and links. We then use these results to give our definition
of Khovanov homology for embedded graphs.
In \cite{kauff} Kauffman introduced a method for producing topological invariants of graphs
embedded in $\s$. The idea is to associate a collection of knots and links to a graph $G$ so that
this family is an invariant under the expanded Reidemeister moves defined by Kauffman and
reported here in figure (\ref{fig21}).
\\
\\
He defined in his work
an ambient isotopy for non-rigid (topological) vertices. (Physically, the rigid vertex
concept corresponds to
a network of rigid disks each with (four) flexible tubes or strings emanating from it.)
Kauffman proved that piecewise linear ambient isotopies of embedded graphs in $\s$
correspond to a sequence of generalized Reidemeister moves for planar diagrams of
the embedded graphs.
\begin{thm}\cite{kauff}
Piecewise linear (PL) ambient isotopy of embedded graphs is generated by the moves of figure (\ref{fig21}), that is, if two embedded graphs
are ambient isotopic, then any two diagrams of them are related by a finite sequence of the moves of figure (\ref{fig21}).
\end{thm}
Let $G$ be a graph embedded in $\s$. The procedure described by Kauffman of how to
associate to $G$ a family of  knots and links prescribes that we should make a local
replacement as in figure (\ref{fig18}) to each vertex in $G$.
Such a replacement at a vertex $v$ connects two edges and isolates all other edges at that vertex, leaving them as free ends. Let $r(G,v)$ denote the link formed by the closed curves formed by this process at a vertex $v$. One retains the link $r(G,v)$, while eliminating all the remaining unknotted arcs.
Define then $T(G)$ to be the family of the links $r(G,v)$ for all possible replacement choices,
$$ T(G)=\cup_{v\in V(G)} r(G,v). $$
For  example see figure (\ref{fig19}).
\\
\begin{thm}\cite{kauff}
Let $G$ be any graph embedded in $\s$, and presented diagrammatically.
Then the family of knots and links $T(G)$,
taken up to ambient isotopy, is a topological invariant of $G$.
\end{thm}

For example, in the figure (\ref{fig19}) the graph $G_2$ is not ambient isotopic to the graph $G_1$,
since $T(G_2)$ contains a non-trivial link.

\subsection{Definition of Khovanov-Kauffman Homology for Embedded Graphs}
Now we are ready to speak about a new concept of Khovanov-Kauffman homology for embedded graphs by using
Khovanov homology for the links (knots) and Kauffman theory of associate a family of links
to an embedded graph $G$, as described above.
\begin{defn}
Let $G$ be an embedded graph with $ T(G)=\{ L_1,L_2,....,L_n\}$ the family of links associated to $G$
by the Kauffman procedure. Let $Kh(L_i)$ be the usual Khovanov homology of the link $L_i$ in this
family. Then the Khovanov-Kauffman homology for the embedded graph $G$ is given by $$ KKh(G)=Kh(L_1) \oplus Kh(L_2)\oplus ....\oplus Kh(L_n)$$
Its graded Euler characteristic is the sum of the graded Euler characteristics of the
Khovanov homology of each link, \ie the sum of the Jones polynomials,
\begin{equation}
\sum_{i,j,k}(-1)^i q^j dim(Kh^{i,j}(L_k))=\sum_k J(L_k).
\end{equation}
\end{defn}
We show some simple explicit examples.
\begin{exa}
In figure (\ref{fig19}) $T(G_1)=\{\bigcirc \bigcirc ,\bigcirc\}$ then for $Kh(\bigcirc)=\mathbb{Q}$
$$ KKh(G_1)=Kh(\bigcirc \bigcirc) \oplus Kh(\bigcirc)$$
Now, from proposition \ref{1234} no.5
$$ KKh(G_1)= Kh(\bigcirc) \otimes Kh(\bigcirc) \oplus Kh(\bigcirc)$$
  $$ KKh(G_1)= \mathbb{Q} \otimes \mathbb{Q} \oplus \mathbb{Q}=\mathbb{Q} \oplus \mathbb{Q}$$
Another example comes from $T(G_2)=\{\includegraphics[width=0.6cm]{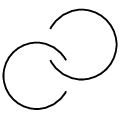}, \bigcirc\}$ then
$$KKh(G_2)=Kh(\includegraphics[width=0.6cm]{hopf.eps}) \oplus Kh(\bigcirc)$$
Since $ Kh^{0,0}(\bigcirc)=\mathbb{Q}$, and from \cite{pt1}

$$Kh(\includegraphics[width=0.6cm]{hopf.eps})= \begin{tabular}{|l|c|c|c|}\hline
\backslashbox{j}{i}  & -2 & -1 & 0 \\\hline
0 &    &      & $\mathbb{Q}$ \\\hline
-1 &    &      &    \\\hline
-2 &    &      & $\mathbb{Q}$    \\\hline
-3 &    &      &    \\\hline
-4 &$\mathbb{Q}$ &      &    \\\hline
-5 &    &      &    \\\hline
-6 &$\mathbb{Q}$ &      &    \\\hline
\end{tabular}
$$

Then,

$$KKh(G_2)= \begin{tabular}{|l|c|c|c|}\hline
\backslashbox{j}{i}  & -2 & -1 & 0 \\\hline
0 &    &      & $\mathbb{Q}\oplus  \mathbb{Q}$ \\\hline
-1 &    &      &    \\\hline
-2 &    &      & $\mathbb{Q}$    \\\hline
-3 &    &      &    \\\hline
-4 &$\mathbb{Q}$ &      &    \\\hline
-5 &    &      &    \\\hline
-6 &$\mathbb{Q}$ &      &    \\\hline
\end{tabular}
$$
\end{exa}

\bibliographystyle{amsplain}

\end{document}